\documentclass[12pt]{amsart}

\usepackage{a4,graphicx,mathrsfs,amssymb, tikz, cancel, color}
\usepackage[colorlinks=true,citecolor=blue]{hyperref}
\usepackage{enumerate}
\usepackage{pstricks-add}
\usepackage{caption}

\newcommand{\C}{{\mathbb C}}
\newcommand{\Z}{{\mathbb Z}}

\newcommand{\N}{{\mathbb N}}
\newcommand{\R}{{\mathbb R}}

\newcommand{\supp}{\mathop{\textrm{supp}}}

\renewcommand{\epsilon}{\varepsilon}

\renewcommand{\tilde}[1]{\widetilde{#1}}
\newtheorem{theorem}{Theorem}[section]

\newtheorem{corollary}[theorem]{Corollary}
\newtheorem{proposition}[theorem]{Proposition}
\newtheorem{remark}[theorem]{Remark}
\newtheorem{definition}[theorem]{Definition}
\theoremstyle{remark}

\numberwithin{equation}{section}

\allowdisplaybreaks

\begin{document}

\title[Equivalent norms in Sobolev spaces]{On some equivalent norms in Sobolev spaces on bounded domains and on the boundaries}\thanks{}

\author{Bienvenido Barraza Mart\'inez}
\address{B.\ Barraza Mart\'inez, Universidad del Norte, Departamento de Matem\'aticas y Estad\'istica, Barranquilla, Colombia}
\email{bbarraza@uninorte.edu.co}

\author{Jonathan Gonz\'alez Ospino}
\address{J.\ Gonz\'alez Ospino, Universidad del Norte, Departamento de Matem\'aticas y Estad\'istica, Barranquilla, Colombia}
\email{gjonathan@uninorte.edu.co}

\author{Jairo Hern\'andez Monz\'on}
\address{J.\ Hern\'andez Monz\'on, Universidad del Norte, Departamento de Matem\'aticas y Estad\'istica, Barranquilla, Colombia}
\email{jahernan@uninorte.edu.co}

\renewcommand{\shortauthors}{B. Barraza Mart\'inez et al.}

\begin{abstract}
We consider the equivalence of some norms in Sobolev spaces on bounded domains of  $\R^d$ and also in Sobolev spaces on the boundaries of those domains.
\end{abstract}



\maketitle

\section{Introduction}
\label{intro}
In this notes we will consider  a bounded domain (i.e., a bounded open and connected subset) $\Omega$ of $\R^d$, $d\in \N$, with enough regular boundary $\partial\Omega$ (this \emph{regularity} will be made precise later). We will present a relativ general result about the equivalence of norms in the scalar Sobolev space $W^{k,p}(\Omega)$ with $k\in\N$ and $1\leq p<\infty$. The main result follows strongly the proof of Theorem 7.1 in \cite{Necas}. For a domain $\Omega\subset\R^d$ (bounded or unbounded) the usual Sobolev space $W^{m,p}(\Omega)$ for $m\in\N$ and $1\leq p\leq\infty$, is the subspace of $L^p(\Omega)$ consisting of all complex fuctions $u\in L^p(\Omega)$ such that its distributional (weak) derivatives $\partial^\alpha u$, with $\alpha\in\N_0^d$ and $\vert\alpha\vert\leq m$, belong to $L^p(\Omega)$. A standard norm in $W^{m,p}(\Omega)$ for $1\leq p<\infty$ is given by
\begin{equation}\label{Eq0}
\Vert u\Vert_{m,p}:=\Vert u\Vert_{m,p,\Omega}:=\bigg( \sum\limits_{\vert\alpha\vert\leq m}\int_{\Omega}\vert\partial^\alpha u(x)\vert^p\,dx\bigg)^{1/p} \qquad (u\in W^{m,p}(\Omega)).
\end{equation}
The Sobolev space $W^{m,\infty}(\Omega)$ is usually endowed with the norm
\begin{equation}\label{Eq0-b}
\Vert u\Vert_{m,\infty} := \Vert u\Vert_{m,\infty,\Omega} := \max\limits_{\vert\alpha\vert\leq m}\Vert \partial^\alpha u\Vert_{L^\infty(\Omega)} \qquad (u\in W^{m,\infty}(\Omega)).
\end{equation}
Now, for $m\in\R$, $m>0$, $m\notin\Z$, and $1\leq p<\infty$, the Sobolev space $W^{m,p}(\Omega)$ (also called \textit{Sobolev-Slobodetskii spaces}) is the subspace of $W^{[m],p}(\Omega)$, where $[m]$ denotes the integer part of $m$, of functions $u$ such that for $\alpha\in \N_0^d$ with $\vert\alpha\vert=[m]$ the following holds
\begin{equation}\label{Eq0-c}
\int_\Omega\int_\Omega \frac{\vert \partial^\alpha u(x)-\partial^\alpha u(y)\vert^p}{\vert x-y\vert^{d+p(m-[m])}}\,dx\,dy < \infty.
\end{equation}
In this case the usual norm in $W^{m,p}(\Omega)$ is given by
\begin{equation}\label{Eq0-d}
\Vert u\Vert_{m,p}:=\Vert u\Vert_{m,p,\Omega}:= \bigg(\Vert u\Vert_{[m],p}^p + \sum\limits_{\vert\alpha\vert=[m]} \int_\Omega\int_\Omega \frac{\vert \partial^\alpha u(x)-\partial^\alpha u(y)\vert^p}{\vert x-y\vert^{d+p(m-[m])}}\,dx\,dy\bigg)^{1/p}
\end{equation}
for $u\in W^{m,p}(\Omega)$.\\
The Sobolev space $W_0^{m,p}(\Omega)$, for $m\geq0$ and $p\geq 1$, is defined as the clousure of $C_c^\infty(\Omega)$ in $W^{m,p}(\Omega)$.\\
Finally, for $m<0$ and $p\geq1$, $W^{m,p}(\Omega)$ is defined as the topological dual of $W_0^{-m,q}(\Omega)$, where $q:=\frac{p}{p-1}$.
\section{About the regularity of the boundaries of domains}\label{Section_regularity_boundary}
In this section we will make precise the concept of continuous and smooth boundary for a bounded domain in $\R^d$. Let $\Omega$ a bounded domain in $\R^d$ with boundary $\partial\Omega$.
\begin{definition}\label{Def-1}
The boundary $\partial\Omega$ is called \emph{continuous} if there exist real numbers $a>0$, $b>0$, a system of local coordinates  $\big\{(x_r^1,\dots,x_r^{d-1},x_r^d)=:(x_r^\prime,x_r^d)\,:\,r=1,\dots,m\big\}$ and continuous functions $a_r:\overline{\Delta_r}\to\R$, $r=1,\dots,m$, where $\Delta_r$ are the open cubes in $\R^{d-1}$ defined by $\Delta_r:=\{(x_r^1,\dots,x_r^{d-1})\,:\,\vert x_r^j\vert<a, j=1,\dots,d-1\}$, such that for each point on the boundary $\partial\Omega$ there is an open neighborhood $V$, such that for some $r\in\{1,\dots,m\}$ the following holds (see Fig. \ref{fig:domain_continuous_boundary})
\begin{align*}
V\cap \partial\Omega & = \big\{(x_r',a_r(x_r'))\,:\,x_r'\in \Delta_r\big\},\\
V\cap \Omega & = \big\{(x_r',x_r^d)\,:\, x_r'\in \Delta_r,\ a_r(x_r')< x_r^d < a_r(x_r')+b \big\},\\
V\cap (\R^d\smallsetminus\overline{\Omega}) & = \big\{(x_r',x_r^d)\,:\, x_r'\in \Delta_r, \ a_r(x_r')-b < x_r^d < a_r(x_r') \big\}.
\end{align*}
\begin{figure}
    \centering
    \includegraphics[scale=0.6]{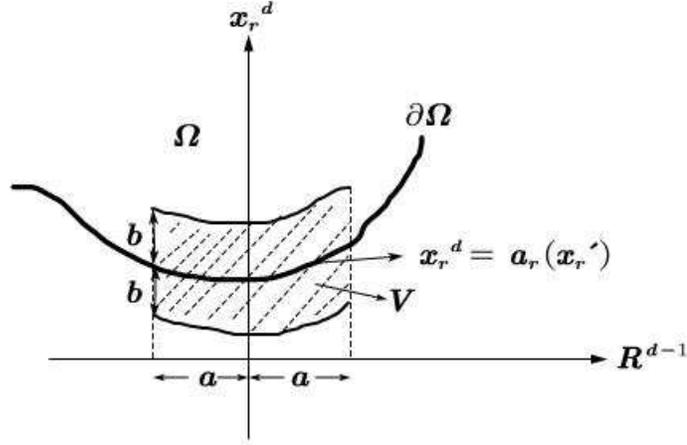}
    \caption{Domain $\Omega$ with continuous boundary $\partial\Omega$.}
    \label{fig:domain_continuous_boundary}
\end{figure}
\end{definition}
\begin{remark}
Note that $\partial\Omega$ is continuos if locally it is the graph of a continuos function defined in a subset of $\R^{d-1}$.
\end{remark}
\begin{definition}
If the functions $a_r$ in the Definition  \ref{Def-1} are Lipschitz continuous it is said that $\Omega$ has a Lipschitz boundary or that the boundary $\partial\Omega$ is lipschitzian. Usually it is said also that the domain $\Omega$ belongs to the class $\mathfrak{N}^{0,1}$.
\end{definition}
\begin{remark}
If $\Omega\in\mathfrak{N}^{0,1}$, then a normal vector exists almost everywhere on $\partial\Omega$ (cf. \cite{Necas}, Lemma 4.2 on pag. 83). In fact, since the functions $a_r$ in Definition \ref{Def-1} are Lipschitz continuous, they are differentiable almost everywhere in their domains. Therefore, for each part $\partial\Omega_r$ of the $\partial\Omega$, represented locally as the graph of the functions $a_r$ for some $r\in\{1,\dots,m\}$, the gradient $\nabla a_r$ exists almost everywhere in $\Delta_r$ and since $\partial\Omega_r$ is a level set of the function $\Delta_r\ni x_r'\mapsto a_r(x_r')-x_r^d$, the vector $(\nabla a_r(x_r'),-1)=\big(\frac{\partial}{\partial x_r^1}a_r(x_r'),\dots,\frac{\partial }{\partial x_r^{d-1}}a_r(x_r'),-1\big)$ defined almost everywhere in $\partial\Omega_r$, is normal to $\partial\Omega_r$ pointing to the exterior of $\Omega$. Then the outer unit normal vector to $\partial\Omega_r$ is given by
\begin{align*}
\nu := \nu(x_r') & := \big(1+\vert\nabla a_r(x_r')\vert^2\big)^{-1/2}\big(\nabla a_r(x_r'),-1\big)\\
& = \bigg(1+\Big(\frac{\partial}{\partial x_r^1}a_r(x_r')\Big)^2 + \cdots +\Big(\frac{\partial}{\partial x_r^{d-1}}a_r(x_r')\Big)^2\bigg)^{-1/2}\\
& \qquad\qquad \cdot\bigg(\frac{\partial}{\partial x_r^1}a_r(x_r'),\dots,\frac{\partial }{\partial x_r^{d-1}}a_r(x_r'),-1\bigg)
\end{align*}
almost everywhere on $\partial\Omega_r$.
\end{remark}
\begin{definition}[Domains of class $\mathfrak{N}^{k,\mu}$]
Let $k\in\N_0$ and $0\leq\mu\leq 1$. It is said that the domain $\Omega$ in Definition \ref{Def-1} belong to the class $\mathfrak{N}^{k,\mu}$ if the functions $a_r$, $r=1,\dots,m$,  given in that definition are of class $C^{k,\mu}(\overline{\Delta_r})$, i.e., if $a_r$ together with its derivatives of order $\leq k$ are H\"older continuous with exponent $\mu$ in $\overline{\Delta_r}$, which means that for each $\alpha\in\N_0^{d-1}$ there is a constant $c$ such that for all $x_r', y_r'\in\overline{\Delta_r}$ the estimate
$$ \vert \partial^\alpha a_r(x_r') - \partial^\alpha a_r(y_r')\vert \leq c\vert x_r' - y_r'\vert^\mu$$
holds.
\end{definition}
\begin{remark}
Note that the case $k=0$ and $\mu=1$ is the Lipschitz case mentioned previously above. In case that $\mu=0$ we also say that $\Omega$ is a domain of  class $C^k$. The notion of continuous boundary given in the first definition of this section corresponds to the case $k=0$ and $\mu=0$.
\end{remark}
\section{Lebesgue and Sobolev spaces on the boundary}\label{Section_Lebesgue_Sobolev_spaces_boundary}
Let $\Omega$ a bounded domain in $\R^d$ with continuous boundary $\partial\Omega$. The notations in this section refer to those given in Definition \ref{Def-1}\footnote{The spaces defined in this section are independent on the local system of coordinates choosen in Definition \ref{Def-1}. The corresponding norms related to each local system of coordinates are all equivalents.}.
\begin{definition}\label{Def-Lp-boundary}
Let $1\leq p\leq \infty$. It is said that a complex function $f$ defined almost everywhere on $\partial\Omega$ (which means that $x_r'\mapsto f(x_r',a_r(x_r'))$ is defined almost everywhere in $\Delta_r$, $r=1,\dots,m$) belongs to the space $L^p(\partial\Omega)$ if the function $x_r'\mapsto f_r(x_r'):=f(x_r',a_r(x_r'))$ belongs to $L^p(\Delta_r)$ for each $r\in\{1,\dots,m\}$. The space $L^p(\partial\Omega)$ is a Banach space with the norm given by
\begin{equation}\label{Norm-Lp_boundary}
\Vert f\Vert_{p,\partial\Omega}:= \Vert f\Vert_{L^p(\partial\Omega)}:= \begin{cases}
\Big(\sum\limits_{r=1}^m \Vert f_r\Vert_{L^p(\Delta_r)}^p\Big)^{1/p}, & \text{if}\ 1\leq p<\infty\\
\max\limits_{r=1,\dots,m}\Vert f_r\Vert_{L^\infty(\Delta_r)}, & \text{if}\ p=\infty.
\end{cases}
\end{equation}
\end{definition}
If $\Omega\in\mathfrak{N}^{0,1}$, the space $L^p(\partial\Omega)$, with $1\leq p<\infty$, can be endowed with another useful norm, equivalent to the norm in \eqref{Norm-Lp_boundary}, which is given in terms of a boundary integral. Next we define the boundary integral for a function in $L^1(\partial\Omega)$. 
\begin{definition}\label{Def_boundary_integral}
Let $\Omega\in \mathfrak{N}^{0,1}$. With the notations of Definition \ref{Def-1}, let
\begin{equation*}
V_r:=\big\{(x_r',x_r^d)\in\R^d\,:\,x_r'\in\Delta_r,\ a_r(x_r')-b<x_r^d<a_r(x_r')+b\big\}, \quad r=1,\dots,m.
\end{equation*}
Now, let $\{\varphi_r\}_{r=1}^m$ a partition of the unity on $\partial\Omega$ subordinate to the cover $\{V_r\}_{r=1}^m$, i.e., for each $r=1,\dots,m$, $\varphi_r\in C_c^\infty(V_r)$, $0\leq \varphi_r\leq 1$, and it holds $\sum\limits_{r=1}^m\varphi_r(x)=1$ for all $x\in\partial\Omega$.\\
For a function $f\in L^1(\partial\Omega)$ we have $f=\sum\limits_{r=1}^m \varphi_r f$ and we define
\begin{equation*}
\int_{\partial\Omega}f\,d\sigma := \sum\limits_{r=1}^m \int_{\Delta_r}f(x_r',a_r(x_r'))\varphi_r(x_r',a_r(x_r'))\sqrt{1+\vert\nabla a_r(x_r')\vert^2}\,dx_r'.
\end{equation*}
\end{definition}
\begin{proposition}
Let $\Omega\in\mathfrak{N}^{0,1}$ and $1\leq p<\infty$. The functional
\begin{equation}\label{Eq_norm_Lp_boundary_integral}
f\mapsto \Big(\displaystyle\int_{\partial\Omega}\vert f\vert^p\,d\sigma \Big)^{1/p}
\end{equation}
is a norm in $L^p(\partial\Omega)$, equivalent to the norm given in \eqref{Norm-Lp_boundary}.
\end{proposition}
\begin{proof}
See \cite{Necas}, Lemma 1.2., pag. 116.
\end{proof}
The following is a standard definiton for the Sobolev spaces on the boundary $\partial\Omega$.
\begin{definition}[Sobolev spaces on the boundary]\label{Def_Sobolev_space_boundary}
Let $k\geq0$, $1\leq p\leq \infty$ and $\Omega\in\mathfrak{N}^{\lceil k\rceil-1,1}$, where $\lceil k\rceil$ is the smallest integer greater than or equal to $k$. The Sobolev space $W^{k,p}(\partial\Omega)$ is the subspace of $L^p(\partial\Omega)$ consisting of all functions $f\in L^p(\partial\Omega)$ such that $f_r\in W^{k,p}(\Delta_r)$ for $r=1,\dots,m$. The space $W^{k,p}(\partial\Omega)$ is endowed with the norm
\begin{equation}\label{Eq_norm_Sobolev_boundary}
\Vert f\Vert_{k,p,\partial\Omega}:= \Vert f\Vert_{W^{k,p}(\partial\Omega)}:= \begin{cases}
\Big(\sum\limits_{r=1}^m \Vert f_r\Vert_{k,p,\Delta_r}^p\Big)^{1/p}, & \text{if}\ 1\leq p<\infty\\
\max\limits_{r=1,\dots,m}\Vert f_r\Vert_{k,\infty,\Delta_r}, & \text{if}\ p=\infty.
\end{cases}
\end{equation}
With this norm $W^{k,p}(\partial\Omega)$ is a Banach space.
\end{definition}
\section{Imbeddings and traces}
\begin{theorem}[Sobolev imbedding theorem]\label{Sobolev-embed-theo}
Let $\Omega\in \mathfrak{N}^{0,1}$, $p\geq 1$ and $kp>d$. Then $W^{k,p}(\Omega)\hookrightarrow C^{0,\mu}(\overline{\Omega})$, where
$$\mu\begin{cases}
= k-\frac{d}{p}, & \text{if}\ k-\frac{d}{p} < 1,\\
< 1, & \text{if}\ k-\frac{d}{p} = 1,\\
= 1, & \text{if}\ k-\frac{d}{p} > 1.
\end{cases}$$
\end{theorem}
\begin{proof}
See Theorem 3.8 in \cite{Necas}, pag. 66.
\end{proof}
\begin{theorem}[A first trace theorem]\label{Trace-theo-1}
Let $\Omega\in\mathfrak{N}^{0,1}$. For $p\leq q\leq \frac{(d-1)p}{d-p}$ if $1\leq p<d$, or for $q\geq 1$ if $p=d$, there exists a continuos linear mapping $Z:W^{1,p}(\Omega)\to L^q(\partial\Omega)$ such that $Zu=u\vert_{\partial\Omega}$ if $u\in C^\infty(\overline{\Omega})$.
\end{theorem}
\begin{proof}
Since $L^{q_2}(\partial\Omega)\hookrightarrow L^{q_1}(\partial\Omega)$ for $1\leq q_1\leq q_2$, the result follows from Theorem 4.2 on page 79 and Theorem 4.6 on page 81 of \cite{Necas}. Cf. also with Theorem 5.36 in \cite{Adams}.
\end{proof}
From Theorems \ref{Sobolev-embed-theo} and \ref{Trace-theo-1} above, we have
\begin{theorem}[A second trace theorem]\label{Trace-theo-2}
Let $\Omega\in\mathfrak{N}^{0,1}$ and $p\geq1$. There is a bounded linear mapping $\gamma_0:W^{1,p}(\Omega)\to L^p(\partial\Omega)$ such that $\gamma_0u=u\vert_{\partial\Omega}$ if 
$u\in C^\infty(\overline{\Omega})$.
\end{theorem}
Usually the mapping $\gamma_0$ is called \textit{trace map of order zero}. The map $\gamma_0:W^{1,p}(\Omega)\to L^p(\partial\Omega)$ is not surjective, but it holds that $\gamma_0(W^{1,p}(\Omega))$ is dense in $L^p(\partial\Omega)$ , whenever $p\geq 1$ and $\Omega\in\mathfrak{N}^{0,1}$ (see Theorem 4.9 on page 82 in \cite{Necas}).\\

A further useful result about traces is the following.
\begin{theorem}[A third trace theorem]\label{Trace-theo-3}
Let $k\in\N$, $p>1$, $\Omega\in\mathfrak{N}^{k-1,1}$ and $u\in W^{k,p}(\Omega)$. If $l\in\N_0$ is such that $l\leq k-1$, then $\gamma_0\partial^\alpha u\in W^{k-l-\frac{1}{p},p}(\partial\Omega)$ for all $\alpha\in\N_0^d$ with $\vert\alpha\vert = l$, and 
\begin{equation}\label{Eq_trace}
\Vert \partial_\nu^l u\Vert_{k-l-\frac{1}{p},p,\partial\Omega}\leq \mathrm{const.}\Vert u\Vert_{k,p,\Omega},
\end{equation}
where 
\begin{equation}\label{Eq_normal derivatives}
\partial_\nu^l u = \sum\limits_{\vert \alpha\vert = l}\frac{l!}{\alpha!}\nu^\alpha\gamma_0\partial^\alpha u,
\end{equation}
where $\nu$ is the outer normal on $\partial\Omega$.
\end{theorem}
\begin{proof}
See \cite{Necas}, Th. 5.5, pag. 95.
\end{proof}
From now on we will omit sometimes the notation $\gamma_0u$ and write simply $u$.
\section{Main results about equivalence of norms in Sobolev spaces on bounded domains}
The following is the main result of this notes.
\begin{theorem}\label{Main theorem}
Let $d, k, l\in\N$, $1\leq p<\infty$, $\Omega\subset \R^d$ a bounded domain with Lipschitz boundary $\partial\Omega$, $\{f_i\}_{i=1}^l$ a set of seminorms in $W^{k,p}(\Omega)$ such that
\begin{itemize}
\item[i)] For each $i=1,\dots,l$, there exists $C_i\geq0$ with $f_i(v)\leq C_i\Vert v\Vert_{k,p}$ for all $v\in W^{k,p}(\Omega)$.
\item[ii)] For $v\in P_{k-1}:=\Big\{\sum\limits_{\vert \alpha\vert\leq k-1}C_\alpha x^\alpha\,:\,C_\alpha\in\R, x\in\R^d\Big\}$ it holds that $\sum\limits_{i=1}^l f_i^p(v) = 0$ implies $v=0$.
\end{itemize}
Then,
\begin{equation}\label{Eq1}
u \mapsto \Vert u\Vert_{k,p}^{\prime} := \bigg(\,\sum\limits_{i=1}^l f_i^p(u) + \vert u\vert_{k,p}^p\bigg)^{1/p},
\end{equation}
with
\begin{equation}\label{Eq2}
\vert u\vert_{k,p}:= \bigg( \sum\limits_{\vert\alpha\vert = k}\int_\Omega \vert \partial^\alpha u (x)\vert^p\,dx\bigg)^{1/p},
\end{equation}
is a norm in $W^{k,p}(\Omega)$, equivalent to the standard one given in \eqref{Eq0}.
\end{theorem}
\begin{proof}
It is clear, due to i), that there exists $b>0$ such that $$\Vert u\Vert_{k,p}^{\prime}\leq b \Vert u\Vert_{k,p}\qquad (u\in W^{k,p}(\Omega)).$$ 
Now suppose that there does not exist a constant $a>0$ such that
$$ a\Vert u\Vert_{k,p} \leq \Vert u\Vert_{k,p}^{\prime}$$
for all $u\in W^{k,p}(\Omega)$. Then, for each $n\in \N$, there exists $u_n\in W^{k,p}(\Omega)$ with $\Vert u_n\Vert_{k,p}=1$ and such that
\begin{equation}\label{Eq3}
\frac{1}{n} > \bigg(\,\sum\limits_{i=1}^l f_i^p(u_n) + \vert u_n\vert_{k,p}^p\bigg)^{1/p}.
\end{equation}
Therefore, for each multiindex $\alpha\in\N_0^d$ with $\vert\alpha\vert=k$, we have
\begin{equation}\label{Eq4}
\partial^\alpha u_n \to 0 \quad \text{in}\ L^p(\Omega), \quad \text{when}\ n\to\infty.
\end{equation}
Theorem 6.3 in \cite{Necas} ensures that the identity mapping $Id: W^{1,p}(\Omega)\to L^p(\Omega)$ is compact, which implies that the identity mapping $Id:W^{k,p}(\Omega)\to W^{k-1,p}(\Omega)$ is also compact. Since $(u_n)_{n\in\N}$ is a bounded sequence in $(W^{k,p}(\Omega), \Vert\cdot\Vert_{k,p})$, there exists a subsequence $(u_{n_m})_{m\in\N}$ of $(u_n)_{n\in\N}$, which converges in the space $(W^{k-1,p}(\Omega),\Vert \cdot\Vert_{k-1,p})$. Let $u:=\lim\limits_{m\to\infty}u_{n_m}$, where the limit is taked in $W^{k-1,p}(\Omega)$. We assert that $u\in W^{k,p}(\Omega)$, $\partial^\alpha u = 0$ for all multiindex $\alpha$ with $\vert\alpha\vert=k$ and that $u_{n_m}\to u$ in $W^{k,p}(\Omega)$ when $m\to\infty$. In fact, let $\varphi\in C_c^\infty(\Omega)$ and $\alpha\in\N_0^d$ with $\vert \alpha\vert = k$. Then, there are $\beta\in\N_0^d$ with $\vert\beta\vert = k-1$ and $j\in\{1,\dots,d\}$ such that $\alpha=\beta + e_j$, where $e_j=(\delta_{ij})_{i=1}^d\in\N_0^d$ with $\delta_{ij}=0$ if $i\neq j$ and $\delta_{ij}=1$ if $i=j$. Due to \eqref{Eq4}  it holds that
\begin{align*}
\int_\Omega u\, \partial^\alpha\varphi\,dx & = \int_\Omega u \,\partial^\beta\partial^{e_j}\varphi\,dx = (-1)^{k-1}\int_\Omega \partial^\beta u\, \partial^{e_j}\varphi\,dx \\
& = (-1)^{k-1}\lim\limits_{m\to\infty}\int_\Omega \partial^\beta u_{n_m}\, \partial^{e_j}\varphi\,dx  = (-1)^{k}\lim\limits_{m\to\infty}\int_\Omega \partial^\alpha u_{n_m}\,\varphi\,dx = 0.
\end{align*}
Then, in weak sense, $\partial^\alpha u = 0\in L^p(\Omega)$. So we have $u\in W^{k,p}(\Omega)$ and, by virtue of \eqref{Eq4}, $\partial^\alpha u_{n_m}\to \partial^\alpha u$ in $L^p(\Omega)$ when $m\to\infty$, for all $\alpha\in\N_0^d$ with $\vert\alpha\vert = k$. Therefore, $u_{n_m}\to u$ in $W^{k,p}(\Omega)$ when $m\to\infty$. Since $\partial^\alpha u=0$ for all $\alpha\in\N_0^d$ with $\vert \alpha\vert = k$, we have that $u\in P_{k-1}$ (see \cite{Dupont-Scott}, Theorem 3.2). Now, since the $f_i$, $i=1,\dots,l$ are seminorms in $W^{k,p}(\Omega)$ and due to the assumption i), for each $i=1,\dots,l$, we have
$$\vert f_i(u_{n_m}) - f_i(u)\vert \leq f_i(u_{n_m}-u)\leq C_i\Vert u_{n_m}-u\Vert_{k,p}\xrightarrow[m\to\infty]{}0.$$ 
By virtue of \eqref{Eq3} we obtain
$$ \bigg(\sum\limits_{i=1}^l f_i^p(u_{n_m})\bigg)^{1/p} < \frac{1}{n_m}.$$
Taking limit when $m\to \infty$ in the last inequality it follows
$$ \sum\limits_{i=1}^l f_i^p(u) = 0,$$
which implies $u=0$ because of the assumption ii). But this is a contradiction with the fact that 
$$\Vert u\Vert_{k,p}=\lim\limits_{m\to\infty}\underbrace{\Vert u_{n_m}\Vert_{k,p}}_{=1} = 1.$$
We came to this contradiction due to the assumption that there does not exist a constant $a>0$ such that $ a\Vert u\Vert_{k,p} \leq \Vert u\Vert_{k,p}^{\prime}$ for all $u\in W^{k,p}(\Omega)$. In consequence, this assumption is false and therefore, there exists $a>0$ such that the inequality
$$ a\Vert u\Vert_{k,p} \leq \Vert u\Vert_{k,p}^{\prime} $$
holds for all $u\in W^{k,p}(\Omega)$. With this we end the proof of Theorem \ref{Main theorem}.
\end{proof}
\begin{remark}\label{remark1}
It is clear that the functional $\Vert \cdot\Vert_{k,p}^\prime$ given in Theorem \ref{Main theorem} is a seminorm in $W^{k,p}(\Omega)$. Theorem 3.2 in \cite{Dupont-Scott} implies that it is in fact a norm in $W^{k,p}(\Omega)$, because if $\Vert u\Vert_{k,p}^\prime = 0$ for $u\in W^{k,p}(\Omega) $, then $\partial^\alpha u=0$ for all $\alpha\in \N_0^d$ with $\vert\alpha\vert = k$ and therefore $u\in P_{k-1}$. We would have also $\sum\limits_{i=1}^l f_i^p(u)=0$ and then $u=0$ by assumption ii).
\end{remark}
\begin{corollary}\label{Corol_5.3}
Let $d, k\in\N$, $1\leq p<\infty$, $\Omega\subset \R^d$ a bounded domain with sufficiently regular boundary $\partial\Omega$ (at least $\Omega\in \mathfrak{N}^{\lceil k-\frac{1}{p}\rceil-1,1}$), $\nu$ the outer normal on $\partial\Omega$, $\Gamma\subseteq \partial\Omega$ with $\sigma(\Gamma)\neq 0$, where $\sigma$ is the $(d-1)$-dimensional Lebesgue surface measure. Furthermore, suppose that $\Gamma$ is not contained in a hyperplane of $\R^d$. Then, the functional
\begin{equation}\label{Eq5}
u\mapsto  \bigg( \sum\limits_{i=0}^{k-1}\int_\Gamma \vert \partial_\nu^i u\vert^p\,d\sigma + \sum\limits_{\vert\alpha\vert=k}\int_\Omega \vert \partial^\alpha u\vert^p\,dx\bigg)^{1/p}
\end{equation}
is a norm in $W^{k,p}(\Omega)$, equivalent to the standard one $\Vert\cdot\Vert_{k,p}$.
\end{corollary}
\begin{proof}
For $i=0,...,k-1$, let $f_i$ be defined by
\begin{equation*}
f_i(u):=\bigg(\int_\Gamma \vert \partial_\nu^i u\vert^p\,d\sigma\bigg)^{1/p} \qquad (u\in W^{k,p}(\Omega)).
\end{equation*}
It is easy to see that the functional $f_i$, $i=0,\dots,k-1$, is a seminorm in $W^{k,p}(\Omega)$. Furthermore, due to equation \eqref{Eq_trace} in Trace theorem \ref{Trace-theo-3}, we have that there exists $C_i\geq0$ such that
\begin{equation*}
f_i(v)\leq \bigg(\int_{\partial\Omega} \vert \partial_\nu^i v\vert^p\,d\sigma\bigg)^{1/p}\leq \mathrm{const.}\Vert \partial_\nu^i v\Vert_{k-i-\frac{1}{p},\, p, \partial\Omega}\leq C_i\Vert v\Vert_{k,p,\Omega}
\end{equation*}
for all $v\in W^{k,p}(\Omega)$. \\
On the other side, let $v\in P_{k-1}$, $v(x)=\sum\limits_{\vert\beta\vert\leq k-1}c_\beta x^\beta$ such that $\sum\limits_{i=0}^{k-1}f_i^p(v)=0$. This implies $\partial_\nu^i v=0$ on $\Gamma$ for each $i=0,\dots,k-1$. We recall that $\partial^\alpha x^\beta = \alpha! \binom{\beta}{\alpha}x^{\beta-\alpha}$ with $\binom{\beta}{\alpha}:=\frac{\beta!}{\alpha!(\beta-\alpha)!}$ if $\alpha\leq\beta$, $\binom{\beta}{\alpha}=0$ otherwise (see \cite{Schleinkofer}, pag. 18). Now, if $\vert\alpha\vert =\vert \beta\vert$ and $\alpha\neq\beta$, then there exists $j\in\{1,\dots,d\}$ such that $\alpha_j>\beta_j$, otherwise $\alpha_i\leq\beta_i$ for all $i\in\{1,\dots,d\}$ and since $\alpha\neq\beta$, we would have $\alpha_l<\beta_l$ for some $l\in\{1,\dots,d\}$ and then $\vert\alpha\vert<\vert\beta\vert$ which contradicts $\vert\alpha\vert=\vert\beta\vert$. Therefore, if $\vert\alpha\vert=\vert\beta\vert$ we have $\partial^\alpha x^\beta = 0$ if $\alpha\neq\beta$ and $\partial^\alpha x^\beta = \alpha!$ if $\alpha=\beta$. We recall also from \eqref{Eq_normal derivatives} that 
$$ \partial_\nu^i v = \sum\limits_{\vert\alpha\vert=i}\frac{i!}{\alpha!}\partial^\alpha v\,\nu^\alpha.$$
Then
\begin{align*}
0 & = \partial_\nu^{k-1}v = \sum\limits_{\vert\alpha\vert = k-1}\frac{(k-1)!}{\alpha!}\partial^\alpha v\,\nu^\alpha = \sum\limits_{\vert\alpha\vert = k-1}\frac{(k-1)!}{\alpha!}\nu^\alpha\sum\limits_{\vert\beta\vert\leq k-1}c_\beta\partial^\alpha x^\beta\\
& = \sum\limits_{\vert\alpha\vert = k-1}\frac{(k-1)!}{\alpha!}\nu^\alpha\sum\limits_{\vert\beta\vert = k-1}c_\beta\partial^\alpha x^\beta = \sum\limits_{\vert\alpha\vert = k-1}\frac{(k-1)!}{\alpha!}\nu^\alpha c_\alpha \alpha!\\
& = \sum\limits_{\vert\alpha\vert = k-1} (k-1)! c_\alpha \nu^\alpha.
\end{align*}
Since $\Gamma$ is not contained in a hiperplane, the powers $\nu^\alpha$ are linear independent and then $c_\alpha =0$ for all $\alpha\in\N_0^d$ with $\vert\alpha\vert=k-1$. Therefore $v(x)=\sum\limits_{\vert\beta\vert \leq k-2} c_\beta x^\beta$. Similarly $\partial_\nu^{k-2} v=0$ implies $c_\beta = 0$ for all $\beta\in\N_0^d$ with $\vert \beta\vert = k-2$. In this form we obtain that $c_\beta=0$ for all $\beta\in \N_0^d$ with $\vert\beta\vert\leq k-1$, i.e., $v=0$.\\
In consequence, we have proved that the functionals $f_i$, $i=0,\dots,k-1$, satisfy the assumptions of Theorem \ref{Main theorem} and we conclude that
\begin{align*}
u & \mapsto \bigg(\,\sum\limits_{i=0}^{k-1} f_i^p(u) + \vert u\vert_{k,p}^p\bigg)^{1/p} = \bigg( \sum\limits_{i=0}^{k-1}\int_\Gamma \vert \partial_\nu^i u\vert^p\,d\sigma + \sum\limits_{\vert\alpha\vert=k}\int_\Omega \vert \partial^\alpha u\vert^p\,dx\bigg)^{1/p}
\end{align*}
is a norm in $W^{k,p}(\Omega)$, equivalent to the norm $\Vert\cdot\Vert_{k,p}$.
\end{proof}
\begin{theorem}[Generalized Poincar\'e inequality]\label{Theo_general_Poincare}
Let $\Omega\subset\R^d$ be open, bounded and connected with Lipschitz boundary $\partial\Omega$ (i.e. $\Omega\in \mathfrak{N}^{0,1}$). Moreover, let $1<p<\infty$ and let $M\subset W^{1,p}(\Omega)$ be nonempty, closed and convex. Then the following assertions are equivalent for every $u_0\in M$:
\begin{enumerate}
\item There exists a constant $C_0<\infty$ such that for all $\xi\in\R$,
$$u_0+\xi\in M \quad \Longrightarrow \quad \vert\xi\vert\leq C_0.$$
\item There exists a constant $C<\infty$ with
$$ \Vert u\Vert_{L^p(\Omega)}\leq C\big(\Vert \nabla u\Vert_{L^p(\Omega)}+1\big)\qquad (u\in M).$$
\end{enumerate}
If $M$ in addition, is a \emph{cone with apex 0}, i.e. if 
$$u\in M, r\geq0 \quad \Longrightarrow\quad ru\in M,$$
then the inequality in the assertion (2) can be replaced with
$$\Vert u\Vert_{L^p(\Omega)}\leq C \Vert \nabla u\Vert_{L^p(\Omega)}\qquad (u\in M).$$
\end{theorem}
\begin{proof}
See \cite[8.16 on pag. 242.]{Alt}
\end{proof}
\begin{corollary}\label{Corol_Poincare_generalized}
Let $\Omega\subset\R^d$ be open, bounded and connected with Lipschitz boundary $\partial\Omega$. Moreover let $1<p<\infty$, $\Gamma\subseteq \partial\Omega$ with $\sigma(\Gamma)\neq 0$, where $\sigma$ is the $(d-1)$-dimensional Lebesgue surface measure, and
$$W_\Gamma^{1,p}(\Omega):=\{u\in W^{1,p}(\Omega)\,:\,u=0 \text{\ on\ } \Gamma\}.$$
Then there exists a constant $C<\infty$ such that
$$ \Vert u\Vert_{L^p(\Omega)}\leq C \Vert \nabla u\Vert_{L^p(\Omega)}\qquad (u\in W_\Gamma^{1,p}(\Omega)).$$
\end{corollary}
\begin{proof}
Let $M$ be defined by
$$M:=\Big\{u\in W^{1,p}(\Omega)\,:\, \int_\Gamma u\,d\sigma = 0\Big\}.$$
Then $M\subset W^{1,p}(\Omega)$ is nonempty because $0\in M$, closed because $u\mapsto  \int_\Gamma u\,d\sigma : W^{1,p}(\Omega)\to \C$ is continuous, and convex because of the linearity of this functional. Now let $u_0\in M$ and take $C_0:=0$. For all $\xi\in\R$ we have
$$u_0+\xi \in M \quad \Longrightarrow \quad 0 = \int_\Gamma (u_0+\xi)\,d\sigma = \int_\Gamma
u_0\,d\sigma + \xi\sigma(\Gamma) = \xi\sigma(\Gamma).$$
Then $\xi=0$ which implies $\vert\xi\vert\leq C_0$. Since $M$ is a cone with appex 0, we have in virtue of Theorem \ref{Theo_general_Poincare} that
$$\Vert u\Vert_{L^p(\Omega)}\leq C \Vert \nabla u\Vert_{L^p(\Omega)}\qquad (u\in M).$$
Let now $u\in W^{1,p}(\Omega)$ and define $\tilde{u}:= u - \dfrac{1}{\sigma(\Gamma)}\displaystyle\int_\Gamma u\,d\sigma$. Then $\tilde{u}\in M$ and we have
$$\Vert \tilde{u}\Vert_{L^p(\Omega)}\leq C \Vert \nabla \tilde{u}\Vert_{L^p(\Omega)}.$$
Then
$$\Vert u\Vert_{L^p(\Omega)} - \frac{\mu(\Omega)^{1/p}}{\sigma(\Gamma)}\Big\vert\int_\Gamma u\,d\sigma \Big\vert \leq \Big\Vert u - \dfrac{1}{\sigma(\Gamma)}\displaystyle\int_\Gamma u\,d\sigma\Big\Vert_{L^p(\Omega)}\leq C \Vert \nabla u\Vert_{L^p(\Omega)},$$
where $\mu(\Omega)$ is the $d$-dimensional Lebesgue measure of $\Omega$. Therefore
$$ \Vert u\Vert_{L^p(\Omega)} \leq C \Vert \nabla u\Vert_{L^p(\Omega)} + \frac{\mu(\Omega)^{1/p}}{\sigma(\Gamma)}\Big\vert\int_\Gamma u\,d\sigma \Big\vert $$
for all $u\in W^{1,p}(\Omega)$. In particular if $u\in W_\Gamma^{1,p}(\Omega)$ we have
$$ \Vert u\Vert_{L^p(\Omega)} \leq C \Vert \nabla u\Vert_{L^p(\Omega)}$$
because in this case $\displaystyle\int_\Gamma u\,d\sigma=0$.
\end{proof}
\begin{remark}
In virtue of Corollary \ref{Corol_Poincare_generalized} we have that the functional 
$u\mapsto \Vert \nabla u\Vert_{L^p(\Omega)}$ is a norm in $W_\Gamma^{1,p}(\Omega)$, which is equivalent to the norm $\Vert\cdot\Vert_{k,p,\Omega}$ (compare with Corollary \ref{Corol_5.3}).
\end{remark}
\section{About equivalence of norms in Sobolev spaces on the boundary}
\begin{theorem}
Let $1\leq p<\infty$ and $\Omega$ a bounded domain in $\R^2$ with sufficiently regular boundary $\partial\Omega$ ( at least $\Omega\in \mathfrak{N}^{2,1}$). Then, the functional
\begin{equation}\label{Eq_equiv_norm_boundary_1}
u\mapsto \big( \Vert u\Vert_{1,p,\partial\Omega}^p + \Vert \partial_\tau^2u\Vert_{p,\partial\Omega}^p\big)^{1/p},
\end{equation}
where $\tau$ is the unit tangential vector on $\partial\Omega$, is a norm in $W^{2,p}(\partial\Omega)$, equivalent to the standard norm $\Vert\cdot\Vert_{2,p,\partial\Omega}$ given in \eqref{Eq_norm_Sobolev_boundary}.
\end{theorem}
\begin{proof}
let $u\in W^{2,p}(\partial\Omega)$. With the notations of Sections \ref{Section_regularity_boundary} and \ref{Section_Lebesgue_Sobolev_spaces_boundary} we have
\begin{align}
\begin{split}\label{Eq_6-2}
\Vert u\Vert_{2,p,\partial\Omega}^p & = \sum\limits_{r=1}^m \Vert u_r\Vert_{2,p,\Delta_r}^p\\
& = \sum\limits_{r=1}^m \sum\limits_{j=0}^2\Big\Vert \Big(\frac{d}{dx_r^1}\Big)^j u_r\Big\Vert_{L^p(\Delta_r)}^p\\
& = \sum\limits_{r=1}^m \sum\limits_{j=0}^1\Big\Vert \Big(\frac{d}{dx_r^1}\Big)^j u_r\Big\Vert_{L^p(\Delta_r)}^p + \sum\limits_{r=1}^m \Big\Vert \Big(\frac{d}{dx_r^1}\Big)^2 u_r\Big\Vert_{L^p(\Delta_r)}^p\\
& = \sum\limits_{r=1}^m \Vert u_r\Vert_{1,p,\Delta_r}^p + \sum\limits_{r=1}^m \Big\Vert \Big(\frac{d}{dx_r^1}\Big)^2 u_r\Big\Vert_{L^p(\Delta_r)}^p\\
& = \Vert u\Vert_{1,p,\partial\Omega}^p + \sum\limits_{r=1}^m \Big\Vert \Big(\frac{d}{dx_r^1}\Big)^2 u_r\Big\Vert_{L^p(\Delta_r)}^p,
\end{split}
\end{align}
where $[x_r^1\mapsto u_r(x_r^1):=u(x_r^1,a_r(x_r^1))]\in W^{2,p}(\Delta_r)=W^{2,p}((-a,a))$, with $a_r\in C^{1,1}(\overline{\Delta_r})=C^{1,1}([-a,a])$, $r=1,\dots,m$. Note that $\Delta_r=(-a,a)$ for all $r=1,\dots,m$.\\
Now, fix $r\in\{1,\dots,m\}$. Taking in account that $\overline{\Delta_r}\ni x_r^1\mapsto (x_r^1,a_r(x_r^1))$ is a parametrization of a part of $\partial\Omega$, we have that $(1,a_r'(x_r^1))$ is a tangent vector to $\partial\Omega$ on that part. Set $\theta_r(x_r^1):=\vert (1,a_r'(x_r^1))\vert = \sqrt{1+(a_r'(x_r^1))^2}$. Furthermore, the weak (or distributional) derivative $\Big(\dfrac{d}{dx_r^1}\Big)^2 u_r$ is almost everywhere equal to the corresponding usual classical derivative in $\Delta_r$ (see Theorem 2.2. in \cite{Necas}, pag. 55). Then, it holds almost everywhere in $\Delta_r$ that
\begin{align}
\begin{split}\label{Eq_6-3}
\Big(\frac{d}{dx_r^1}\Big)^2 & u_r(x_r^1)\\
& = \frac{d}{dx_r^1}\big[\nabla u (x_r^1,a_r(x_r^1))\cdot (1,a_r'(x_r^1))\big]\\
& = \frac{d}{dx_r^1}\big[\theta_r(x_r^1)\partial_\tau u (x_r^1,a_r(x_r^1))\big]\\
& = \theta_r'(x_r^1)\partial_\tau u (x_r^1,a_r(x_r^1)) + \theta_r(x_r^1)\nabla(\partial_\tau u)(x_r^1,a_r(x_r^1))\cdot (1,a_r'(x_r^1))\\
& =  \theta_r'(x_r^1)\partial_\tau u (x_r^1,a_r(x_r^1)) + \theta_r(x_r^1)^2\partial_\tau^2 u(x_r^1,a_r(x_r^1))\\
& = \frac{\theta_r'(x_r^1)}{\theta_r(x_r^1)}\frac{d}{dx_r^1}u(x_r^1,a_r(x_r^1)) + \theta_r(x_r^1)^2\partial_\tau^2 u(x_r^1,a_r(x_r^1))\\
& = \frac{\theta_r'(x_r^1)}{\theta_r(x_r^1)}\frac{d}{dx_r^1}u_r(x_r^1) + \theta_r(x_r^1)^2\partial_\tau^2 u(x_r^1,a_r(x_r^1)).
\end{split}
\end{align}
Since $a_r\in C^{2,1}(\overline{\Delta_r})$, there are constants $c_r^1$, $c_r^2$ and $c_r^3$ such that
\begin{align}
\begin{split}\label{Eq_constantes}
\bigg\vert\dfrac{\theta_r'(x_r^1)}{\theta_r(x_r^1)}\bigg\vert\leq c_r^1,\\
0 < c_r^2 \leq \theta_r(x_r^1)^2\leq c_r^3.
\end{split}
\end{align} 
Then, due to \eqref{Eq_6-3} we have with $[\partial_\tau^2u]_r(x_r^1):= \partial_\tau^2 u(x_r^1,a_r(x_r^1))$ that
\begin{align*}
\Big\Vert \Big(\frac{d}{dx_r^1}\Big)^2 u_r\Big\Vert_{L^p(\Delta_r)} & \leq \Big\Vert \frac{\theta_r'}{\theta_r}\frac{d}{dx_r^1}u_r\Big\Vert_{L^p(\Delta_r)} + \big\Vert \theta_r^2[\partial_\tau^2 u]_r\big\Vert_{L^p(\Delta_r)}\\
& \leq c_r^1\Big\Vert \frac{d}{dx_r^1}u_r\Big\Vert_{L^p(\Delta_r)} + c_r^3\big\Vert [\partial_\tau^2 u]_r\big\Vert_{L^p(\Delta_r)}.
\end{align*}
Therefore
\begin{equation*}
\Big\Vert \Big(\frac{d}{dx_r^1}\Big)^2 u_r\Big\Vert_{L^p(\Delta_r)}^p  \leq 2^p(c_r^1)^p\Big\Vert \frac{d}{dx_r^1}u_r\Big\Vert_{L^p(\Delta_r)}^p + 2^p(c_r^3)^p\big\Vert [\partial_\tau^2 u]_r\big\Vert_{L^p(\Delta_r)}^p.
\end{equation*}
From \eqref{Eq_6-2}, with $c_1:=\max\limits_{r=1,\dots,m}2^p(c_r^1)^p$ and $c_3:=\max\limits_{r=1,\dots,m}2^p(c_r^3)^p$, we obtain
\begin{align*}
\Vert u\Vert_{2,p,\partial\Omega}^p & \leq  \Vert u\Vert_{1,p,\partial\Omega}^p + c_1\sum\limits_{r=1}^m \Big\Vert \frac{d}{dx_r^1}u_r\Big\Vert_{L^p(\Delta_r)}^p + c_3\sum\limits_{r=1}^m \big\Vert [\partial_\tau^2 u]_r\big\Vert_{L^p(\Delta_r)}^p\\
& \leq \Vert u\Vert_{1,p,\partial\Omega}^p + c_1\sum\limits_{r=1}^m \Vert u_r\Vert_{1,p,\Delta_r}^p + c_3\sum\limits_{r=1}^m \big\Vert [\partial_\tau^2 u]_r\big\Vert_{L^p(\Delta_r)}^p\\
& = (1+c_1)\Vert u\Vert_{1,p,\partial\Omega}^p + c_3\Vert \partial_\tau^2u\Vert_{p,\partial\Omega}^p\\
& \leq \max\{1+c_1, c_3\}\big(\Vert u\Vert_{1,p,\partial\Omega}^p + \Vert \partial_\tau^2u\Vert_{p,\partial\Omega}^p \big).
\end{align*}
With $\tilde{c}_1:=(\max\{1+c_1, c_3\})^{1/p}$, it holds
\begin{equation}\label{Eq_6-4}
\Vert u\Vert_{2,p,\partial\Omega} \leq \tilde{c}_1\big(\Vert u\Vert_{1,p,\partial\Omega}^p + \Vert \partial_\tau^2u\Vert_{p,\partial\Omega}^p \big)^{1/p}.
\end{equation}
On the other side, again from \eqref{Eq_6-3}, we have
\begin{align*}
c_r^2\big\Vert[\partial_\tau^2 u]_r&\big\Vert_{L^p(\Delta_r)}\\
 & \leq \big\Vert\theta_r^2[\partial_\tau^2 u]_r\big\Vert_{L^p(\Delta_r)} \\
 & \leq \Big\Vert \Big(\frac{d}{dx_r^1}\Big)^2 u_r\Big\Vert_{L^p(\Delta_r)} + c_r^1\Big\Vert \frac{d}{dx_r^1} u_r\Big\Vert_{L^p(\Delta_r)},
\end{align*}
which implies
\begin{align*}
(c_r^2)^p\big\Vert[\partial_\tau^2 u]_r&\big\Vert_{L^p(\Delta_r)}^p\\
& \leq 2^p\Big\Vert \Big(\frac{d}{dx_r^1}\Big)^2 u_r\Big\Vert_{L^p(\Delta_r)}^p + 2^p(c_r^1)^p\Big\Vert \frac{d}{dx_r^1} u_r\Big\Vert_{L^p(\Delta_r)}^p,
\end{align*}
i.e.,
\begin{align*}
\big\Vert[\partial_\tau^2 u]_r&\big\Vert_{L^p(\Delta_r)}^p\\
& \leq \max\bigg\{ \frac{2^p}{(c_r^2)^p} ,\frac{2^p(c_r^1)^p}{(c_r^2)^p}\bigg \}\bigg(\Big\Vert \Big(\frac{d}{dx_r^1}\Big)^2 u_r\Big\Vert_{L^p(\Delta_r)}^p + \Big\Vert \frac{d}{dx_r^1} u_r\Big\Vert_{L^p(\Delta_r)}^p\bigg)\\
& \leq \max\bigg\{ \frac{2^p}{(c_r^2)^p} ,\frac{2^p(c_r^1)^p}{(c_r^2)^p}\bigg \}\Vert u_r\Vert_{2,p,\Delta_r}^p.
\end{align*}
Therefore, with $c_2:=\max\limits_{r=1,\dots,m}\max\bigg\{ \frac{2^p}{(c_r^2)^p} ,\frac{2^p(c_r^1)^p}{(c_r^2)^p}\bigg \}$, we have
\begin{align*}
\Vert u\Vert_{1,p,\partial\Omega}^p & + \Vert \partial_\tau^2u\Vert_{p,\partial\Omega}^p\\ 
& = \Vert u\Vert_{1,p,\partial\Omega}^p + \sum\limits_{r=1}^m\big\Vert[\partial_\tau^2 u]_r\big\Vert_{L^p(\Delta_r)}^p\\
 & \leq \Vert u\Vert_{1,p,\partial\Omega}^p + c_2\sum\limits_{r=1}^m \Vert u_r\Vert_{2,p,\Delta_r}^p\\
 & = \Vert u\Vert_{1,p,\partial\Omega}^p + c_2\Vert u\Vert_{2,p,\partial\Omega}^p\\
 & \leq (1+c_2)\Vert u\Vert_{2,p,\partial\Omega}^p.
\end{align*}
With $\tilde{c}_2:=(1+c_2)^{1/p}$ it holds
\begin{equation}\label{Eq_6-6}
\Big(\Vert u\Vert_{1,p,\partial\Omega}^p  + \Vert \partial_\tau^2u\Vert_{p,\partial\Omega}^p\Big)^{1/p} \leq \tilde{c}_2\Vert u\Vert_{2,p,\partial\Omega}.
\end{equation}
From \eqref{Eq_6-4} and \eqref{Eq_6-6} follows the result.
\end{proof}
\section{Some useful interpolation type estimates for traces on Sobolev spaces}
In this section we use as norm in $L^p(\partial\Omega)$, for a domain $\Omega\subset\R^d$, the norm given in \eqref{Eq_norm_Lp_boundary_integral}.
\begin{figure}
    \centering
    \includegraphics[scale=0.45]{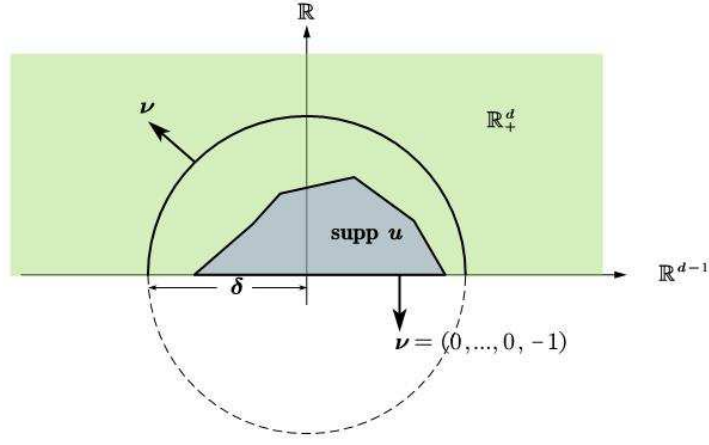}
    \caption{Domain $\R_+^d:=\{x=(x',x_d)\in\R^{d-1}\times\R\,:\,x_d>0\}$ and $\supp u$.}
    \label{fig:domain_R_mas_supp_u}
\end{figure}
\begin{proposition}\label{Prop_interp_estimates_1}
Let $\R_+^d:=\big\{x=(x',x_d)\in\R^{d-1}\times\R\,:\,x_d>0\big\}$  and $u\in C_c^1(\overline{\R_+^d})$ with $1\leq p<\infty$. The following estimate holds
\begin{equation}\label{Eq_estimate_interp_1}
\Vert u\Vert_{p,\partial\R_+^d} \leq p^{1/p}\,\Vert u\Vert_{p,\R_+^d}^{\frac{p-1}{p}}\Vert u\Vert_{1,p,\R_+^d}^{\frac{1}{p}}.
\end{equation}
\end{proposition}
\begin{proof}
Let $\delta>0$ be such that $\supp u\subset B(0;\delta)\cap \overline{\R_+^d}$ and $\nu=(\nu_1,\dots,\nu_d)$ the outer normal on $\partial[B(0;\delta)\cap {\R_+^d}]$ (see Fig. \ref{fig:domain_R_mas_supp_u}). Then, by virtue of Gau\ss\  theorem of divergence and H\"older inequality it holds
\begin{align*}
\int_{\partial\R_+^d}\vert u\vert^p\,d\sigma & = \int_{\R^{d-1}}\vert u\vert^p\,dx' = - \int_{\R^{d-1}}\vert u\vert^p(-1)\,dx' = -\, \int_{\partial[B(0;\delta)\cap {\R_+^d}]}\vert u\vert^p\nu_d\,d\sigma\\
& = -\, \int_{B(0;\delta)\cap {\R_+^d}}\partial_{x_d}(\vert u\vert^p)\,dx = - \int_{\R_+^d}\partial_{x_d}(u\overline{u})^{p/2}\,dx\\
& = -\,\frac{p}{2}\int_{\R_+^d}(u\overline{u})^{\frac{p-2}{2}}[(\partial_{x_d}u)\overline{u} + u\partial_{x_d}\overline{u}]\,dx\\
& = -\, p\int_{\R_+^d}\vert u\vert^{p-2}\mathrm{Re}[(\partial_{x_d}u)\overline{u}]\,dx = -\,p\,\mathrm{Re} \int_{\R_+^d}\vert u\vert^{p-2}(\partial_{x_d}u)\overline{u}\,dx\\
& \leq p\,\int_{\R_+^d}\vert u\vert^{p-2}\vert \partial_{x_d}u\vert \vert u\vert\,dx = p\,\int_{\R_+^d}\vert u\vert^{p-1}\vert \partial_{x_d}u\vert\,dx\\
& \leq p\,\Big(\int_{\R_+^d}\vert u\vert^p\,dx\Big)^{\frac{p-1}{p}}\Big(\int_{\R_+^d}\vert \partial_{x_d}u\vert^p\,dx\Big)^{\frac{1}{p}}\\
& \leq p\,\Vert u\Vert_{p,\R_+^d}^{p-1}\Vert u\Vert_{1,p,\R_+^d},
\end{align*} 
which implies \eqref{Eq_estimate_interp_1}.
\end{proof}
\begin{proposition}\label{Prop_interp_estimates_2}
Let $\Omega$ a bounded domain in $\R^d$, $\Omega\in\mathfrak{N}^{1,1}$, $1\leq p<\infty$ and $u\in C^1(\overline{\Omega})$. Then, it holds
\begin{equation}\label{Eq_estimate_interp_2}
\Vert u\Vert_{p,\partial\Omega} \leq c_p\,\Vert u\Vert_{p,\Omega}^{\frac{p-1}{p}}\Vert u\Vert_{1,p,\Omega}^{\frac{1}{p}},
\end{equation}
with $c_p$ a constant, which depends on $p$, but not on $u$.
\end{proposition}
\begin{proof}
With the notations of Section \ref{Section_regularity_boundary}, let 
\begin{equation*}
V_r:=\big\{(x_r',x_r^d)\in\R^d\,:\,x_r'\in\Delta_r,\ a_r(x_r')-b<x_r^d<a_r(x_r')+b\big\}, \quad r=1,\dots,m.
\end{equation*}
Furthermore let $V_0\subset\subset\Omega$, $V_0$ abierto, such that $\overline{\Omega}\subset\bigcup\limits_{r=0}^mV_r$. Choose  a $C^\infty$ partition of the unity $\{\varphi_r\}_{r=0}^m$ on $\overline{\Omega}$ subordinate to the cover $\{V_r\}_{r=0}^m$, i.e., $\varphi_r\in C_c^\infty(V_r)$, $0\leq \varphi_r\leq 1$ and $\supp \varphi_r\subset V_r$ for $r=0,1,\dots,m$. Moreover $\sum\limits_{r=0}^m\varphi_r(x) = 1$ for all $x\in\overline{\Omega}$. Then, $u=\sum\limits_{r=0}^m u\varphi_r$ on $\overline{\Omega}$ with $u\varphi_r\in C_c^1(\overline{V_r^+})$ for $r=0,1,\dots,m$, where $V_r^+:=V_r\cap\Omega$.\\

Now, for $r=1,\dots,m$, we will consider the transformation of coordinates $V_r\ni (x_r',x_r^d)\longleftrightarrow (y_r',y_r^d)$ given by
\begin{equation*}
y_r=\Phi_r(x_r):=\begin{cases}
y_r^i := x_r^i, & i=1,\dots,d-1,\\
y_r^d := x_r^d - a_r(x_r'), &
\end{cases}
\end{equation*}
with inverse transformation
\begin{equation*}
x_r=\Phi_r^{-1}(y_r):=\begin{cases}
x_r^i = y_r^i, & i=1,\dots,d-1,\\
x_r^d = y_r^d + a_r(y_r'). &
\end{cases}
\end{equation*}
Let 
$$w^r(y_r) := u(y_r',y_r^d + a_r(y_r'))\varphi_r(y_r',y_r^d + a_r(y_r')), \quad r=1,\dots,m,$$
for $y_r\in\Phi_r(\overline{V_r^+})$. We have $w^r\in C_c^1(\Phi_r(\overline{V_r^+}))$ and therefore, extending by zero outside of $\Phi_r(\overline{V_r^+})$, $w^r\in C_c^1(\overline{\R_+^d})$. Then, it follows (with several constants $c_r$, $c^1$, $c_p^2$, etc., which can depend on $p$, but not on $u$)
\begin{align*}
\Vert u\Vert_{p,\partial\Omega}^p & = \Big\Vert \sum\limits_{r=0}^m u\varphi_r\Big\Vert_{p,\partial\Omega}^p \leq m\sum\limits_{r=1}^m\Vert u\varphi_r\Vert_{p,\partial\Omega}^p = m\sum\limits_{r=1}^m\int_{\partial\Omega}\vert u\varphi_r\vert^p\,d\sigma\\
& = m\sum\limits_{r=1}^m \int_{\Delta_r}\vert u(x_r',a_r(x_r'))\varphi_r(x_r',a_r(x_r'))\vert^p\sqrt{1+\vert\nabla a_r(x_r')\vert^2}\,dx_r'\\
& \leq m\sum\limits_{r=1}^m c_r\int_{\R^{d-1}}\vert u(y_r',a_r(y_r'))\varphi_r(y_r',a_r(y_r'))\vert^p\,dy_r'\\
& \leq c^1\sum\limits_{r=1}^m \int_{\R^{d-1}} \vert w^r(y_r',0)\vert^p\,dy_r'\\
& \leq c^1p\sum\limits_{r=1}^m \Vert w^r\Vert_{p,\R_+^d}^{p-1}\Vert w^r\Vert_{1,p,\R_+^d} \qquad\qquad (\text{Prop. \ref{Prop_interp_estimates_1}})\\
& = c_p^2\sum\limits_{r=1}^m \Vert w^r\Vert_{p,\Phi_r({V_r^+})}^{p-1}\Vert w^r\Vert_{1,p,\Phi_r({V_r^+})}\\
& \leq c_p^3\sum\limits_{r=1}^m \Vert u\varphi_r\Vert_{p,V_r^+}^{p-1}\Vert u\varphi_r\Vert_{1,p,V_r^+} \qquad\qquad (\text{Th. 4.1 in \cite{Wloka}, p. 80}).
\end{align*}
Now,  
$$\Vert u\varphi_r\Vert_{p,V_r^+}\leq \Vert u\Vert_{p,\Omega}$$
and 
$$ \Vert u\varphi_r\Vert_{1,p,V_r^+}^p = \Vert u\varphi_r\Vert_{p,V_r^+}^p + \Vert \nabla(u\varphi_r)\Vert_{p,V_r^+}^p \leq c_r^4\Vert u\Vert_{1,p,\Omega}^p. $$
Therefore,
$$ \Vert u\Vert_{p,\partial\Omega}^p \leq c_p^3\sum\limits_{r=1}^m (c_r^4)^{1/p}\Vert u\Vert_{p,\Omega}^{p-1}\Vert u\Vert_{1,p,\Omega}.$$
From this follows \eqref{Eq_estimate_interp_2} with $c_p:=\Big(c_p^3\sum\limits_{r=1}^m (c_r^{4})^{1/p}\Big)^{1/p}$.
\end{proof}
\begin{proposition}\label{Prop_interp_estimates_3}
Let $\Omega$ a bounded domain in $\R^d$, $\Omega\in\mathfrak{N}^{1,1}$, $1\leq p<\infty$ and $u\in W^{1,p}(\Omega)$. Then, it holds
\begin{equation}\label{Eq_estimate_interp_3}
\Vert u\Vert_{p,\partial\Omega} \leq c_p\,\Vert u\Vert_{p,\Omega}^{\frac{p-1}{p}}\Vert u\Vert_{1,p,\Omega}^{\frac{1}{p}},
\end{equation}
with $c_p$ a constant, which depends on $p$, but not on $u$. 
\end{proposition}
\begin{proof}
Let $(u_n)_{n\in\N}$ a sequence of functions of $C^\infty(\overline{\Omega})$ such that $u_n\to u$ in $W^{1,p}(\Omega)$ whenever $n\to\infty$. Due to trace theorem \ref{Trace-theo-2} we have
that there exists a constant $c$,  such that
$$\Vert u_n-u\Vert_{p,\partial\Omega}\leq c \Vert u_n-u\Vert_{1,p,\Omega}\qquad (n\in\N).$$
Then, $\Vert u_n-u\Vert_{p,\partial\Omega}\to 0$ when $n\to\infty$.\\

Now, from Proposition \ref{Prop_interp_estimates_2} it follows that
$$ \Vert u_n\Vert_{p,\partial\Omega}\leq c_p\,\Vert u_n\Vert_{p,\Omega}^{\frac{p-1}{p}}\Vert u_n\Vert_{1,p,\Omega}^{\frac{1}{p}}\qquad (n\in\N).$$
Making $n\to\infty$ we obtain \eqref{Eq_estimate_interp_3}.
\end{proof}
From Proposition \ref{Prop_interp_estimates_3} follow also the following estimates.
\begin{proposition}\label{Lemma_02}
  Let $\Omega$ a bounded domain in $\R^d$, $\Omega\in\mathfrak{N}^{2,1}$, $1\leq p<\infty$ and $u\in W^{2,p}(\Omega)$.  The following estimate holds:
  \begin{equation}\label{Ineq_02}
   \Vert \partial_\nu u\Vert_{p,\partial\Omega} \leq \widetilde{c}_p \Vert u\Vert_{1,p,\Omega}^{\frac{p-1}{p}}\Vert u\Vert_{2,p,\Omega}^{\frac{1}{p}}
  \end{equation}
with $\widetilde{c}_p$ being a positive constant independent of $u$, where $\nu$ is the outer normal on $\partial\Omega$.
\end{proposition}
\begin{proof}
 For $u\in C^2(\overline{\Omega})$, due to Proposition \ref{Prop_interp_estimates_3}, the following estimates hold:
 \begin{align*}
  \Vert \partial_\nu u\Vert_{p,\partial\Omega} & = \Big\Vert \sum\limits_{j=1}^d \nu_j\partial_ju\Big\Vert_{p,\partial\Omega} \leq \sum\limits_{j=1}^d \big\Vert  \nu_j\partial_ju \big\Vert_{p,\partial\Omega} \leq \max\limits_{\partial\Omega}\vert \nu\vert\sum\limits_{j=1}^d \big\Vert \partial_ju \big\Vert_{p,\partial\Omega}\\
  & \leq c_p\max\limits_{\partial\Omega}\vert \nu\vert\sum\limits_{j=1}^d\Vert \partial_ju\Vert_{p,\Omega}^{\frac{p-1}{p}}\Vert \partial_ju\Vert_{1,p,\Omega}^{\frac{1}{p}} \leq \widetilde{c}_p\Vert u\Vert_{1,p,\Omega}^{\frac{p-1}{p}}\Vert u\Vert_{2,p,\Omega}^{\frac{1}{p}},
 \end{align*}
where $\widetilde{c}_p:= d\,c_p\max\limits_{\partial\Omega}\vert \nu\vert$, with $c_p$ the constant of Proposition \ref{Prop_interp_estimates_3}. Then, \eqref{Ineq_02} is also true for $u\in W^{2,p}(\Omega)$ due to the density of $C^2(\overline{\Omega})$ in $W^{2,p}(\Omega)$.
\end{proof}

\begin{proposition}\label{Lemma_03}
Let $\Omega$ a bounded domain in $\R^d$, $\Omega\in\mathfrak{N}^{3,1}$, $1\leq p<\infty$ and $u\in W^{3,p}(\Omega)$. The following estimate holds:
 \begin{equation}\label{Ineq_03}
   \Vert \Delta u\Vert_{p,\partial\Omega} \leq \hat{c}_p \Vert u\Vert_{2,p,\Omega}^{\frac{p-1}{p}}\Vert u\Vert_{3,p,\Omega}^{\frac{1}{p}}
  \end{equation}
  with $\hat{c}_p$ being a positive constant independent of $u$.
\end{proposition}
\begin{proof}
Let $u\in C^3(\overline{\Omega})$. Then, similarly as the proof of Proposition \ref{Lemma_02} and using again Propostion \ref{Prop_interp_estimates_3} (or Proposition \ref{Prop_interp_estimates_2}), we have
\begin{align*}
\Vert \Delta u\Vert_{p,\partial\Omega} & = \Big\Vert \sum\limits_{j=1}^d \partial_j^2u\Big\Vert_{p,\partial\Omega} \leq \sum\limits_{j=1}^d \big\Vert \partial_j^2u \big\Vert_{p,\partial\Omega}\\
& \leq c_p\sum\limits_{j=1}^d \big\Vert \partial_j^2u \big\Vert_{p,\Omega}^{\frac{p-1}{p}}\big\Vert \partial_j^2u \big\Vert_{1,p,\Omega}^{\frac{1}{p}} \leq d\,c_p\Vert u\Vert_{2,p,\Omega}^{\frac{p-1}{p}}\Vert u\Vert_{3,p,\Omega}^{\frac{1}{p}}.
\end{align*}
Due to the density of $C^3(\overline{\Omega})$ in $W^{3,p}(\Omega)$ we obtain that the estimate \eqref{Ineq_03} holds also for $u\in W^{3,p}(\Omega)$, with $\hat{c}_p:=d\,c_p$, where $c_p$ is the constant of Proposition \ref{Prop_interp_estimates_3}.
\end{proof}
\begin{proposition}\label{Lemma_04}
 Let $\Omega$ a bounded domain in $\R^d$, $\Omega\in\mathfrak{N}^{4,1}$, $1\leq p<\infty$ and $u\in W^{4,p}(\Omega)$. The following estimate holds:
 \begin{equation}\label{Ineq_04}
   \Vert \partial_\nu \Delta u\Vert_{p,\partial\Omega} \leq \widetilde{c}_p \Vert u\Vert_{3,p,\Omega}^{\frac{p-1}{p}}\Vert u\Vert_{4,p,\Omega}^{\frac{1}{p}}
  \end{equation}
  with $\widetilde{c}_p$ being a positive constant independent of $u$.
\end{proposition}
\begin{proof}
 Similarly as the proof of Proposition \ref{Lemma_02} we obtain that for $u\in C^4(\overline{\Omega})$ the estimate \eqref{Ineq_04} holds. Then, due to the density of $C^4(\overline{\Omega})$ in $W^{4,p}(\Omega)$, the estimate \eqref{Ineq_04} holds also for $u\in W^{4,p}(\Omega)$, with the constant $\widetilde{c}_p$ being the same of Proposition \ref{Lemma_02}.
\end{proof}

%

\end{document}